\documentclass[a4paper,11pt]{article}

\RequirePackage[T1]{fontenc}
\RequirePackage{fix-cm}

\usepackage{graphicx, url}
\usepackage{amsmath,amsfonts,amsthm}
\usepackage{booktabs}
\usepackage{tikz, pgfplots}
\pgfplotsset{compat=1.17}
\usepgfplotslibrary{groupplots}
\usetikzlibrary{shapes.geometric}
\usetikzlibrary{backgrounds}
\usetikzlibrary{intersections}
\usetikzlibrary{spy}
\newcommand{\ddt}{\partial_t}
\newcommand{\ddx}{\partial_x}
\newcommand{\R}{\mathbb{R}}

\newcommand{\eps}{\varepsilon}
\begin{document}

\title{Numerical relaxation limit and outgoing edges in a central scheme for networked conservation laws}

\author{Niklas Kolbe$^{1,}$%
\footnote{Corresponding author; e-mail: \tt{kolbe@igpm.rwth-aachen.de}}}

\date{
  \small
  $^1$Institute of Geometry and Practical Mathematics, RWTH Aachen University,\\ Templergraben 55, 52062 Aachen, Germany\\
   \smallskip
   \today
}
\maketitle                   

\begin{abstract}
A recently introduced scheme for networked conservation laws is analyzed in various experiments. The scheme makes use of a novel relaxation approach that governs the coupling conditions of the network and does not require a solution of the Riemann problem at the nodes. We numerically compare the dynamics of the solution obtained by the scheme to solutions obtained using a classical coupling condition. In particular, we investigate the case of two outgoing edges in the Lighthill--Whitham--Richards model of traffic flow and in the Buckley--Leverett model of two phase flow. Moreover, we numerically study the asymptotic preserving property of the scheme by comparing it to its preliminary form before the relaxation limit in a 1-to-1 network.
\end{abstract}

\section{Networks of scalar conservation laws and the central scheme}

In this work a network refers to a directed graph consisting of edges and nodes. On each of its edges a scalar hyperbolic conservation law is imposed. Due to the finite speed of propagation we restrict the discussion to networks containing only a single coupling node that is placed at position $x=0$. The dynamics of the state variable on the adjacent edge $k$ are governed by
\begin{equation}\label{eq:scalarconservationnetwork}
  \ddt u^k + \ddx f_k(u^k) =0 \quad \text{in }\mathcal{E}_k\times(0,\infty),\quad k \in \delta^\mp,
\end{equation}
where the edge is either incoming and parameterized by $\mathcal{E}_k=(-1,0)$ if $k\in \delta^-=\{1,\dots,N^- \}$ or outgoing and parameterized by $\mathcal{E}_k=(0, 1)$ if $k\in\delta^+=\{N^-+1, \dots,N\}$. The flux functions $f_1,\dots,f_{N}:$ $\R \rightarrow \R$ are smooth, but not-necessarily convex or concave. The set of all edges is denoted by $\delta^{\mp}=\delta^- \cup \delta^+$. In addition, we impose no-flux conditions on the left boundaries and homogeneous Neumann conditions on the right boundaries, i.e.,
\begin{equation}\label{eq:bcs}
f_k(\cdot) = 0 \quad \text{at } x=-1 \quad \text{for } k\in \delta^-, \qquad \ddx u = 0 \quad \text{at } x=1 \quad \text{for }k\in \delta^+   
\end{equation}
and we assume given initial data denoted on each edge $k\in \delta^\mp$ by $u^{k,0}$. Applications of such network models include gas dynamics in pipelines~\cite{banda2006coupleuler}, vehicular traffic on road networks~\cite{holden1995} and production systems \cite{dapice2010model}. For well-posedness of the network problem \eqref{eq:scalarconservationnetwork}, \eqref{eq:bcs} coupling conditions at the node are required.

We are concerned with the numerical approximation of the network \eqref{eq:scalarconservationnetwork}. Rewriting the conservation laws in terms of the relaxation system from \cite{jin1995relaxschemsystem}, discretizing by an asymptotic-preserving scheme and taking the relaxation limit at the coupling node, a suitable scheme has been introduced in \cite{herty2022centr}. This scheme governs the coupling conditions of the network and is Riemann solver-free. It can be written in conservative form as
\begin{equation}\label{eq:conservativeformnet}
  u_j^{k,n+1} = u_j^{k,n} - \frac{\Delta t }{\Delta x} \left( F_{j+1/2}^{k,n} - F_{j-1/2}^{k,n}\right),
\end{equation}
where $u_j^{k,n}$ is an approximate average of $u^k$ over the cell $I_j=[(j-1/2)\Delta x, (j+1/2)\Delta x]$ at time instance $t= n\,\Delta t$ using the time- and space increments $\Delta t>0$ and $\Delta x = \frac 1 m$. In our numerical experiments we choose the time increment according to the CFL condition $\Delta t = \text{CFL}\, \frac{\Delta x}{\max_k \lambda_k}$. The subindices in \eqref{eq:conservativeformnet} can be taken $j\in \{ -m,\dots, -1\}$ if $k\in \delta^-$ or $j\in \{ 0, 1, \dots, m\}$ if $k\in \delta^+$. The numerical fluxes are given by
\begin{equation}\label{eq:netfluxes}
  F_{j-1/2}^{k,n} =
  \begin{cases}
    \frac 1 2 \, (f_k(u_{j}^{k,n}) + f_k(u_{j-1}^{k,n}))  - \frac {\lambda_k} 2 (u_j^{k,n}-u_{j-1}^{k,n}) - {\mathcal{S}}^{k,n}_{j-1/2} & \text{if }j \neq 0,\\[5pt]
    \frac 1 2 \, ( v_R^{k,n} + f_k(u_{-1}^{k,n}))  - \frac {\lambda_k} 2 (u_R^{k,n}-u_{-1}^{k,n}) &\text{if }j=0 \text{ and } k \in \delta^-, \\[5pt]
    \frac 1 2 \, ( f_k(u_{0}^{k,n}) + v_L^{k,n})  - \frac {\lambda_k} 2 (u_0^{k,n}-u_{L}^{k,n}) &\text{if }j=0 \text{ and } k\in \delta^+. 
  \end{cases}
\end{equation}
The relaxation speeds $\lambda_k>0$ are chosen such that the subcharacteristic condition $-\lambda_k\leq f_k^\prime(u^k) \leq \lambda_k$ holds at all edges.
Taking $\mathcal{S}^{k,n} =0$ yields a first order scheme, whereas the term can also account for a piecewise linear second order approximation, see \cite{herty2022centr} for details. The coupling data $u_R^{k,n}$, $u_L^{k,n}$, $v_R^{k,n}$, $v_L^{k,n}$ play a key role in the numerical fluxes at the coupling node. They depend on the cell averages next to the coupling node and can be generally obtained from the solution of two linear systems. These systems depend on the structure of the network and allow implementing priority rules between the edges. In case of a 1-to-1 network (i.e. $N^-=N^+=1$) the systems read
  \begin{align}
    \label{eq:couplingdataonetoone}
    \begin{pmatrix}
      \lambda_1 & -\lambda_2 \\[5pt]
      -\lambda_1^2 & \lambda_2^2
    \end{pmatrix}
    \begin{pmatrix}
      u_R^1 - u_{-1} ^1 \\[5pt]
      u_L^2 - u_0^2
    \end{pmatrix}
    &=
    \begin{pmatrix}
      v_{-1}^1 - v_{0}^2 \\[5pt]
      \lambda_1^2 \, u_{-1}^1 - \lambda_2^2 \, u_0^2
    \end{pmatrix}, \\
    \begin{pmatrix}
      -1 & 1 \\[5pt]
      \lambda_1 & \lambda_2
    \end{pmatrix}
    \begin{pmatrix}
      v_R^1 - f_1(u_{-1} ^1) \\[5pt]
      v_L^2 - f_2(u_0^2)
    \end{pmatrix}
    &=
    \begin{pmatrix}
      v_{-1}^1 - v_{0}^2 \\[5pt]
      \lambda_1^2 \, u_{-1}^1 - \lambda_2^2 \, u_0^2
    \end{pmatrix}.
  \end{align}
  Mass conservation of this scheme for any network structure has been verified and an implementation is available from \cite{kolbe2022implem}. 
  
  \section{Numerical relaxation limit}
The study of \emph{asymptotic preserving} schemes has been of high interest in the last decades, see \cite{hu2017asymppreserschem}. Such schemes have the property to preserve the asymptotic transition from an underlying micro-model to the macro one at the discrete level. The scheme introduced above has been derived from a relaxation approach, which takes the role of the micro-model in this context. In the approach instead of the conservation law the relaxation system
  \begin{subequations}
    \begin{align}
      \ddt u + \ddx v &= 0 && \text{in }\R \times (0, \infty), \label{eq:relu} \\
      \ddt v + \lambda^2 \ddx u &= \frac{1}{\eps} (f(u) - v) && \text{in }\R \times (0, \infty)\label{eq:relv}
    \end{align}
  \end{subequations}
  is imposed on the edges of the network and the relaxation limit $\eps\rightarrow 0$ has been taken. Following the steps in \cite{herty2022centr} we derive a scheme for the relaxation network in case $\eps>0$. It includes an additional evolution formula for the cell averages of the auxiliary variable $v$ that reads
  \begin{equation}\label{eq:conservativeformv}
    v_j^{k,n+1} = v_j^{k,n} - \frac{\Delta t }{\Delta x} \left( G_{j+1/2}^{k,n} - G_{j-1/2}^{k,n}\right) + \frac{\Delta t}{\eps} \left(f_k(u^{k,n+1}_j) - v^{k,n+1}_j\right).
  \end{equation}
 To obtain the correct relaxation limit, the source term of the system is implicitly considered. As the state variable $u$ can be updated before the auxiliary variable $v$ it is not necessary to solve a nonlinear system to evaluate \eqref{eq:conservativeformv}. The numerical fluxes are given by
  \begin{equation}\label{eq:vfluxes}
    G_{j-1/2}^{k,n} =
  \begin{cases}
    \frac {\lambda_k^2} 2 \, (u_{j}^{k,n} + u_{j-1}^{k,n})  - \frac {\lambda_k} 2 (v_j^{k,n}-v_{j-1}^{k,n}) & \text{if }j \neq 0,\\[5pt]
    \frac {\lambda_k^2} 2 \, ( u_R^{k,n} + u_{-1}^{k,n})  - \frac {\lambda_k} 2 (v_R^{k,n}-v_{-1}^{k,n}) &\text{if }j=0 \text{ and } k \in \delta^-, \\[5pt]
    \frac {\lambda_k^2} 2 \, ( u_{0}^{k,n} + u_L^{k,n})  - \frac {\lambda_k} 2 (v_0^{k,n}-v_{L}^{k,n}) &\text{if }j=0 \text{ and } k\in \delta^+. 
  \end{cases}
\end{equation}
To account for the state variable $u$ in the relaxation network we impose an update formula given by \eqref{eq:conservativeformnet} and a modification of \eqref{eq:netfluxes} where all terms of the form $f_k(u^{k,n}_j)$ are replaced by $v^{k,n}_j$. The same substitution is considered in the computation of the coupling data according to \eqref{eq:couplingdataonetoone}. Although a second order version of the approximation \eqref{eq:conservativeformv}, \eqref{eq:vfluxes} can be constructed in analogy to \eqref{eq:conservativeformnet}, \eqref{eq:netfluxes} we consider a first order scheme for simplicity and set $\mathcal{S}^{k,n}=0$.

For our numerical study we consider the Lighthill--Whitham--Richards (LWR) model of traffic flow, see \cite{garavello2006traffflownetwor}, on a 1-to-1 network. The numerical experiment is adapted from \cite{chiarello2020} and assumes the flux functions
\[
f_1(u) = 2 \, u (1-2u), \quad f_2(u) = u(1-u)
\]
and the piecewise constant initial data $u^{1,0}\equiv \frac 1 4$ left and $u^{2,0}\equiv \frac 1 2$ right from the coupling node. We employ both schemes, the \emph{limit scheme} introduced in Section 1 and the \emph{relaxation scheme} for the relaxation network, to numerically solve the problem. In the relaxation scheme the initial data of the auxiliary variable is chosen as $v^{1,0}= f_1(u^{1,0})$ and $v^{2,0}= f_2(u^{2,0})$. As we focus on the role of relaxation we choose a fine grid employing $m=2000$ cells on both edges. We moreover fix the relaxation speed $\lambda_1=\lambda_2=\lambda=2$ and the Courant number $\text{CFL}=0.9$.

\begin{figure}
  \centering
  \begin{tikzpicture}
  \begin{groupplot}[
        group style={group size=4 by 1,
            horizontal sep = .3 cm, 
            vertical sep = .2 cm,
            xticklabels at=edge bottom,
            yticklabels at=edge left}, 
          width = .31 \linewidth,
          height = .2 \linewidth,
          xmin=-1,xmax=1,ymin=0, ymax=0.6,
          every tick label/.append style={font=\scriptsize},
          title style={font=\scriptsize},
          every axis plot/.append style={
            only marks, mark=+, mark size=.5 pt, line width=0.15
          }
          ]
        \nextgroupplot[title={$t=0$}]
        \addplot [] table [x index=0, y index=1] {input/lwr11_limit_1.dat};
        \nextgroupplot[title={$t=0.5$}]
        \addplot [] table [x index=0, y index=1] {input/lwr11_limit_2.dat};
        \nextgroupplot[title={$t=1$}]
        \addplot [] table [x index=0, y index=1] {input/lwr11_limit_3.dat};
        \nextgroupplot[title={$t=1.5$}]
        \addplot [] table [x index=0, y index=1] {input/lwr11_limit_4.dat};
  \end{groupplot}
\end{tikzpicture}
  \caption{Numerical solution of the first order central scheme applied to the LWR model on a 1-to-1 network. The vehicle densities on the incoming ($x<0$) and the outgoing ($x>0$) road are shown in four time instances. Congestion on the outgoing road causes the formation of a shock that propagates along both roads as time evolves. Computations employed $2000$ mesh cells on both edges/roads and $\text{CFL}=0.9$. }\label{fig:LWR11dynamics} 
\end{figure}
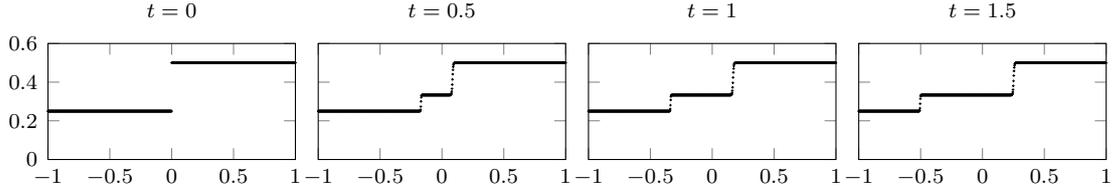

\begin{figure}
  \centering
  \begin{tikzpicture}[spy using outlines={rectangle, red, magnification=5,
    size=1.5cm, connect spies}]
  \begin{groupplot}[
        group style={group size=2 by 1,
            horizontal sep = .35 cm, 
            vertical sep = .2 cm,
            xticklabels at=edge bottom},
          width = .535 \linewidth,
          height = .3 \linewidth,
          xmin=-1,xmax=1,
          every tick label/.append style={font=\scriptsize},
          title style={font=\small},
          legend to name=leg:relaxcompare,
          legend columns=-1,
          legend style={font=\scriptsize, fill=none, /tikz/every even column/.append style={column sep=.25cm}}
          ]
        \nextgroupplot[title={$u$} ,ymin=0.2, ymax=0.6]
        \addplot [cyan, semithick] table [x index=0, y index=1] {input/lwr11_relaxation_1.dat};
        \addplot [magenta, semithick] table [x index=0, y index=1] {input/lwr11_relaxation_2.dat};
        \addplot [lime, semithick] table [x index=0, y index=1] {input/lwr11_relaxation_3.dat};
        \addplot [olive, semithick] table [x index=0, y index=1] {input/lwr11_relaxation_4.dat};
        \addplot [black, semithick] table [x index=0, y index=1] {input/lwr11_limit_4.dat};
        \coordinate (spypoint) at (axis cs:-0.53,0.26);
        \coordinate (magpoint) at (axis cs:-.7,.49);
        \spy [red] on (spypoint) in node at (magpoint);
        \nextgroupplot[title={$v$}, ymin=0.19, ymax=0.26, yticklabel pos=right]
        \addplot [cyan, semithick] table [x index=0, y index=2] {input/lwr11_relaxation_1.dat};
        \addplot [magenta, semithick] table [x index=0, y index=2] {input/lwr11_relaxation_2.dat};
        \addplot [lime, semithick] table [x index=0, y index=2] {input/lwr11_relaxation_3.dat};
        \addplot [olive, semithick] table [x index=0, y index=2] {input/lwr11_relaxation_4.dat};
        \addplot [black, semithick] table [x index=0, y index=2] {input/lwr11_limit_4.dat};
        \coordinate (spypoint2) at (axis cs:0.23,0.225);
        \coordinate (magpoint2) at (axis cs:.7,.21);
        \spy [red] on (spypoint2) in node at (magpoint2);
        \legend{$\eps = 10^{-1}$, $\eps = 10^{-2}$, $\eps = 10^{-3}$, $\eps = 10^{-4}$, $\eps\rightarrow0$}
  \end{groupplot}
\end{tikzpicture}
\ref{leg:relaxcompare}
  \caption{Comparison of the relaxation scheme for varied $\eps$ to the limit scheme ($\eps\rightarrow 0$) at time instance $T=1.5$ in terms of numerically computed vehicle density (left) and auxiliary variable (right). In case of the limit scheme $f_k(u^k)$ is shown instead of the auxiliary variable. Selected regions of the numerical solutions are magnified. The numerical solution of the relaxation scheme approaches the one of the limit scheme as $\eps$ decreases. Computations employed $2000$ mesh cells on both edges left and right from the coupling node ($x=0$) and $\text{CFL}=0.9$.}\label{fig:LWR11comparison} 
\end{figure}
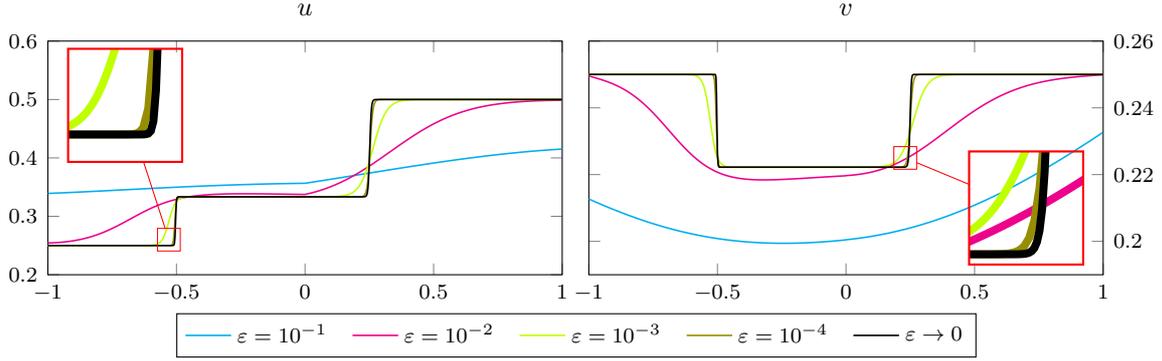

\begin{table}
  \centering
  \scriptsize
  \begin{tabular}{c @{\hspace{4em}} c c c c c c}
    \toprule
    $\eps$ & $10^{-1}$ & $10^{-2}$ & $10^{-3}$ & $10^{-4}$ & $10^{-5}$  & $10^{-6}$ \\ \midrule
    $L^1$ error & $1.414 \times 10^{-1}$ & $4.964 \times 10^{-2}$ & $7.983 \times 10^{-3}$ & $8.182 \times 10^{-4}$& $8.186 \times 10^{-5}$ & $8.187 \times 10^{-6}$\\
    EOC & & $0.45$ & $0.79$ & $0.99$ & $1.00$ &$1.00$  \\ \bottomrule
  \end{tabular}
  \caption{$L^1$ errors comparing the vehicle density obtained by the relaxation scheme for varied $\eps$ to the vehicle density by the limit scheme and EOCs. The $L^1$ error decreases and the EOC indicates a first order convergence.}\label{tab:eps}
\end{table}

Fig.~\ref{fig:LWR11dynamics} shows the numerical solution computed by the limit scheme. A shock wave of vehicle density $\frac{1}{3}$ emerges at the coupling node and propagates to both sides because of congestion on the outgoing road/edge. Due to the different flux functions on the edges the absolute shock velocity on the incoming edge is larger than the one on the outgoing edge. In Fig.~\ref{fig:LWR11comparison} we compare numerical solutions of the relaxation scheme for varied $\eps$ to the numerical solution of the limit scheme at the final time $T=1.5$. Clearly the traffic density obtained by the relaxation scheme tends to the one by the limit scheme as $\eps$ decreases. Similarly, the auxiliary variable obtained by the relaxation scheme tends to $f_k(u^k_\text{lim})$ as $\eps$ goes to zero with $u^k_\text{lim}$ denoting the traffic density obtained by the limit scheme. Furthermore, we compute the $L^1$ error at $T=1.5$ comparing the vehicle density obtained by the relaxation scheme for various $\eps$ to the vehicle density by the limit scheme. Along with these errors that we denote by $E_k$ for $\eps=10^{-k}$ we compute the experimental order of convergence (EOC) by the formula $\log_{10}(E_{k-1}/E_k)$. The computed errors and EOCs shown in Table \ref{tab:eps} confirm the limit behavior observed in Fig. \ref{fig:LWR11comparison} and indicate a first order convergence with respect to $\eps$.

\section{Outgoing edges}
In this section we consider numerical experiments on 1-to-2 networks (i.e., $N^-=1$, $N^+=2$) using the scheme given by \eqref{eq:conservativeformnet} and \eqref{eq:netfluxes}. In this case the coupling data required in the numerical fluxes are computed from the linear systems
\begin{align}
  \label{eq:couplingdataonetoone}
  \begin{pmatrix}
    \lambda_1 & -\lambda_2 & -\lambda_3\\[5pt]
    -\lambda_1^2 & \lambda_2^2 & \lambda_3^2\\[5pt]
    -\alpha \lambda_1 & -\lambda_2 & 0
  \end{pmatrix}
  \begin{pmatrix}
    u_R^1 - u_{-1} ^1 \\[5pt]
    u_L^2 - u_0^2 \\[5pt]
    u_L^3 - u_0^3 \\[5pt]
  \end{pmatrix}
  &=
    \begin{pmatrix}
      v_{-1}^1 - v_{0}^2 - v_{0}^3\\[5pt]
      \lambda_1^2 \, u_{-1}^1 - \lambda_2^2 \, u_0^2 - \lambda_3^2 \, u_0^3 \\[5pt]
      -\alpha v_{-1}^1+v_0^2
    \end{pmatrix}, \\
    \begin{pmatrix}
      -1 & 1 & 1\\[5pt]
      \lambda_1 & \lambda_2 & \lambda_3 \\[5pt]
      \alpha & -1 & 0
    \end{pmatrix}
    \begin{pmatrix}
      v_R^1 - f_1(u_{-1} ^1) \\[5pt]
      v_L^2 - f_2(u_0^2) \\[5pt]
      v_L^2 - f_2(u_0^2)
    \end{pmatrix}
    &=
    \begin{pmatrix}
      v_{-1}^1 - v_{0}^2 - v_{0}^3\\[5pt]
      \lambda_1^2 \, u_{-1}^1 - \lambda_2^2 \, u_0^2 - \lambda_3^2 \, u_0^3 \\[5pt]
      -\alpha v_{-1}^1+v_0^2
    \end{pmatrix},
  \end{align}
  where the parameter $\alpha\in[0,1]$ determines the flux distribution to the outgoing edges. In more details, $\alpha$ is the rate of flux going into edge 2 in the total flux, consequently the rate of flux going into edge 3 in the total flux is $1-\alpha$.

\subsection{Coupled two-phase flow model} 
In this numerical experiment we impose the Buckley--Leverett equation \cite{buckley1942mechanfluiddisplsands}, which models two-phase flow, on a 1-to-2 network. The state variable $u$ takes the role of the water fraction in a mixture of water and oil, which is governed by a conservation law with the non-convex flux function
\begin{equation}\label{eq:bl}
  f_i(u) = \frac{u^2}{u^2 + a_i (1-u)^2},
\end{equation}
where $a_i<1$ is a constant. Our numerical experiment reproduces a scenario, in which water is pumped into a tube to displace oil and enforce its outflow through a second and a third tube. We set the fluxes of the network edges representing the tubes to \eqref{eq:bl} with parameters $a_1= 0.5$, $a_2=0.1$ and $a_3=0.9$ to account for different tube/edge properties. We assume constant initial data on each edge given by $u^{1, 0}\equiv1$, $u^{2,0}\equiv 0$ and $u^{3,0}\equiv 0$.

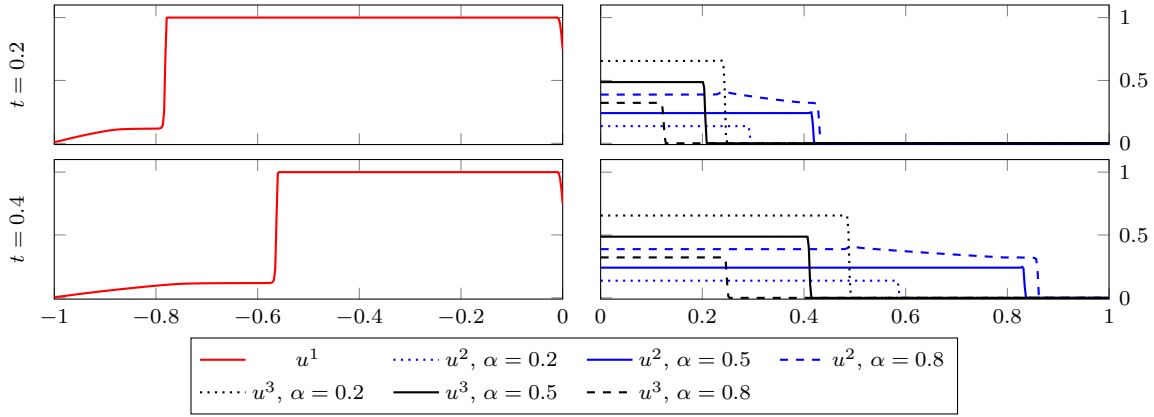
\begin{figure}
  \centering
  \begin{tikzpicture}
  \begin{groupplot}[
        group style={group size=2 by 2,
            horizontal sep = .5 cm, 
            vertical sep = .2 cm,
            xticklabels at=edge bottom,
            yticklabels at=edge right}, 
          width = .53 \linewidth,
          height = .22 \linewidth,
          ymin=-.01, ymax=1.1,
          every tick label/.append style={font=\scriptsize},
          legend to name=leg:congestion,
          legend columns=4,
          legend style={font=\scriptsize, fill=none, /tikz/every even column/.append style={column sep=.3cm}},
          label style={font=\scriptsize}
          ]
        \nextgroupplot[ylabel={$t=0.2$}, xmin=-1, xmax=0]
        \addplot [color=red, thick] table [x index=0, y index=3] {input/BL12_a5_2.dat};
        \nextgroupplot[xmin=0, xmax=1]
        \addplot [color=blue, thick, dotted] table [x index=1, y index=4] {input/BL12_a2_2.dat};
        \addplot [color=blue, thick] table [x index=1, y index=4] {input/BL12_a5_2.dat};
        \addplot [color=blue, thick, dashed] table [x index=1, y index=4] {input/BL12_a8_2.dat};
        \addplot [color=black, thick, dotted] table [x index=1, y index=5] {input/BL12_a2_2.dat};
        \addplot [color=black, thick] table [x index=1, y index=5] {input/BL12_a5_2.dat};
        \addplot [color=black, thick, dashed] table [x index=1, y index=5] {input/BL12_a8_2.dat};
        \nextgroupplot[ylabel={$t=0.4$}, xmin=-1, xmax=0]
        \addplot [color=red, thick] table [x index=0, y index=3] {input/BL12_a5_3.dat};
        \nextgroupplot[xmin=0, xmax=1]
        \addplot [color=red, thick] coordinates {(-5,0)};
        \addplot [color=blue, thick, dotted] table [x index=1, y index=4] {input/BL12_a2_3.dat};
        \addplot [color=blue, thick] table [x index=1, y index=4] {input/BL12_a5_3.dat};
        \addplot [color=blue, thick, dashed] table [x index=1, y index=4] {input/BL12_a8_3.dat};
        \addplot [color=black, thick, dotted] table [x index=1, y index=5] {input/BL12_a2_3.dat};
        \addplot [color=black, thick] table [x index=1, y index=5] {input/BL12_a5_3.dat};
        \addplot [color=black, thick, dashed] table [x index=1, y index=5] {input/BL12_a8_3.dat};
        \legend{$u^1$, $u^{2}\text{, }\alpha=0.2$, $u^{2}\text{, }\alpha=0.5$, $u^{2}\text{, }\alpha=0.8$, $u^{3}\text{, }\alpha=0.2$, $u^3\text{, }\alpha=0.5$, $u^3\text{, }\alpha=0.8$}
      \end{groupplot}
    \end{tikzpicture}
    \ref{leg:congestion}
  \caption{Incoming (left) and outgoing (right) edges/tubes in the 1-to-2 Buckley--Leverett network in two time instances for varied flux distribution parameter $\alpha$. The water fraction in tube 1 is independent of $\alpha$. Smaller $\alpha$ lead to faster water outflow in tube 3. The second order scheme with $\text{CFL}=0.24$ and $400$ grid cells per edge was used for the numerical simulation.}\label{fig:BL12} 
\end{figure}

Fig.~\ref{fig:BL12} shows the time evolution of the numerical solution for varied flux distribution parameter $\alpha$ computed by the second order scheme on $400$ grid cells per edge, relaxation speeds $\lambda_1=\lambda_2=\lambda_3=2.5$ and Courant number $\text{CFL}=0.24$. Shock waves are formed at the coupling node that move along the outgoing edges. Tube 2 allows for a faster outflow of the water. In case of equal flux distribution ($\alpha=0.5$) this causes a faster shock wave in edge 2 but a higher amplitude in tube 3. Reducing the outflow in tube 2 by setting $\alpha=0.2$ increases both the velocity and the amplitude of the shock wave in tube 3 while it reduces the water flow in tube 2. An increased outflow in tube 3 ($\alpha=0.8$ shown) has the inverse effect. In this case the shock wave in tube 2 is followed by a rarefaction wave.

\subsection{Traffic flow}
In the last numerical experiment we impose the LWR model, which we also considered in Section 2, on the edges of a 1-to-2 network. We assume that a larger road that allows for unidirectional traffic flow splits into two smaller ones. The higher capacity of the incoming road is reflected in the flux functions, which we choose
\begin{equation}
  f_1(u)=u \left( 1-\frac{u}{1.2} \right), \quad f_2(u)=f_3(u)=u(1-u).
\end{equation}
A scenario with high vehicle density on the roads that might lead to congestion given by the initial data $u^{1, 0}\equiv0.6$, $u^{2,0}\equiv 0.9$ and $u^{3,0}\equiv 0.4$ is considered. Analogously to the experiments in \cite{herty2022centr}, we compare the numerical solution of the central approach in our scheme to the one obtained by \emph{flow maximization}. On 1-to-2 networks flow maximization leads to the following fluxes at the coupling node
\begin{equation}\label{eq:flowmax}
  F^1_{-1/2} = \min\left\{ d_1(u_{-1}^1), \frac{s_2(u_0^2)}{\alpha}, \frac{s_3(u_0^3)}{1-\alpha}\right\}, \quad F^2_{-1/2} = \alpha \, F^1_{-1/2}, \quad F^3_{-1/2} = (1-\alpha) \, F^1_{-1/2},
\end{equation}
where time indices are neglected and $d_1$, $s_2$ and $s_3$ denote the demand and supply functions corresponding to the flux functions, see~\cite{garavello2006traffflownetwor}. To numerically simulate flow maximization on the network, we replace in the scheme the numerical fluxes at the coupling node by~\eqref{eq:flowmax}.

\begin{figure}
  \centering
  \begin{tikzpicture}
  \begin{groupplot}[
        group style={group size=2 by 2,
            horizontal sep = .5 cm, 
            vertical sep = .2 cm,
            xticklabels at=edge bottom,
            yticklabels at=edge right}, 
          width = .53 \linewidth,
          height = .22 \linewidth,
          ymin=-.01, ymax=1.1,
          every tick label/.append style={font=\scriptsize},
          legend to name=leg:congestion,
          legend columns=5,
          legend style={font=\scriptsize, fill=none, /tikz/every even column/.append style={column sep=.5cm}},
          label style={font=\scriptsize}
          ]
        \nextgroupplot[ylabel={central}, xmin=-1, xmax=0]
        \addplot [color=red, thick, dotted] table [x index=0, y index=3] {input/LWR12_CentralRelaxationLimit_a2_3.dat};
        \addplot [color=red, thick, dashed] table [x index=0, y index=3] {input/LWR12_CentralRelaxationLimit_a4_3.dat};
        \addplot [color=red, thick] table [x index=0, y index=3] {input/LWR12_CentralRelaxationLimit_a8_3.dat};
        \nextgroupplot[xmin=0, xmax=1]
        \addplot [color=blue, thick, dotted] table [x index=1, y index=4] {input/LWR12_CentralRelaxationLimit_a2_3.dat};
        \addplot [color=blue, thick, dashed] table [x index=1, y index=4] {input/LWR12_CentralRelaxationLimit_a4_3.dat};
        \addplot [color=blue, thick] table [x index=1, y index=4] {input/LWR12_CentralRelaxationLimit_a8_3.dat};
        \addplot [color=black, thick, dotted] table [x index=1, y index=5] {input/LWR12_CentralRelaxationLimit_a2_3.dat};
        \addplot [color=black, thick, dashed] table [x index=1, y index=5] {input/LWR12_CentralRelaxationLimit_a4_3.dat};
        \addplot [color=black, thick] table [x index=1, y index=5] {input/LWR12_CentralRelaxationLimit_a8_3.dat};
        \nextgroupplot[ylabel={flow maximization}, xmin=-1, xmax=0]
        \addplot [color=red, thick, dotted] table [x index=0, y index=3] {input/LWR12_TrafficFlowMaximization_a2_3.dat};
        \addplot [color=red, thick, dashed] table [x index=0, y index=3] {input/LWR12_TrafficFlowMaximization_a4_3.dat};
        \addplot [color=red, thick] table [x index=0, y index=3] {input/LWR12_TrafficFlowMaximization_a8_3.dat};
        \nextgroupplot[xmin=0, xmax=1]
        \addplot [color=red, thick, dotted] coordinates {(-5,0)};
        \addplot [color=red, thick, dashed] coordinates {(-5,0)};
        \addplot [color=red, thick]  coordinates {(-5,0)};
        \addplot [color=blue, thick, dotted] table [x index=1, y index=4] {input/LWR12_TrafficFlowMaximization_a2_3.dat};
        \addplot [color=blue, thick, dashed] table [x index=1, y index=4] {input/LWR12_TrafficFlowMaximization_a4_3.dat};
        \addplot [color=blue, thick] table [x index=1, y index=4] {input/LWR12_TrafficFlowMaximization_a8_3.dat};
        \addplot [color=black, thick, dotted] table [x index=1, y index=5] {input/LWR12_TrafficFlowMaximization_a2_3.dat};
        \addplot [color=black, thick, dashed] table [x index=1, y index=5] {input/LWR12_TrafficFlowMaximization_a4_3.dat};
        \addplot [color=black, thick] table [x index=1, y index=5] {input/LWR12_TrafficFlowMaximization_a8_3.dat};
        \legend{$u^1\text{, }\alpha=0.2$, $u^1\text{, }\alpha=0.4$, $u^1\text{, }\alpha=0.8$, $u^{2}\text{, }\alpha=0.2$, $u^{2}\text{, }\alpha=0.4$, $u^{2}\text{, }\alpha=0.8$, $u^{3}\text{, }\alpha=0.2$, $u^3\text{, }\alpha=0.4$, $u^3\text{, }\alpha=0.8$}
      \end{groupplot}
    \end{tikzpicture}
    \ref{leg:congestion}
  \caption{Incoming (left) and outgoing (right) edges in the 1-to-2 LWR network computed using both the central (top) and the flow maximization (bottom) approach for varied flux distribution parameter $\alpha$. Both coupling approaches predict qualitatively similar dynamics. Larger $\alpha$ lead to congestion and a backward moving traffic wave on the incoming road. The second order scheme with $\text{CFL}=0.24$ and $200$ grid cells per edge was used for the numerical simulation.}\label{fig:LWR12} 
\end{figure}
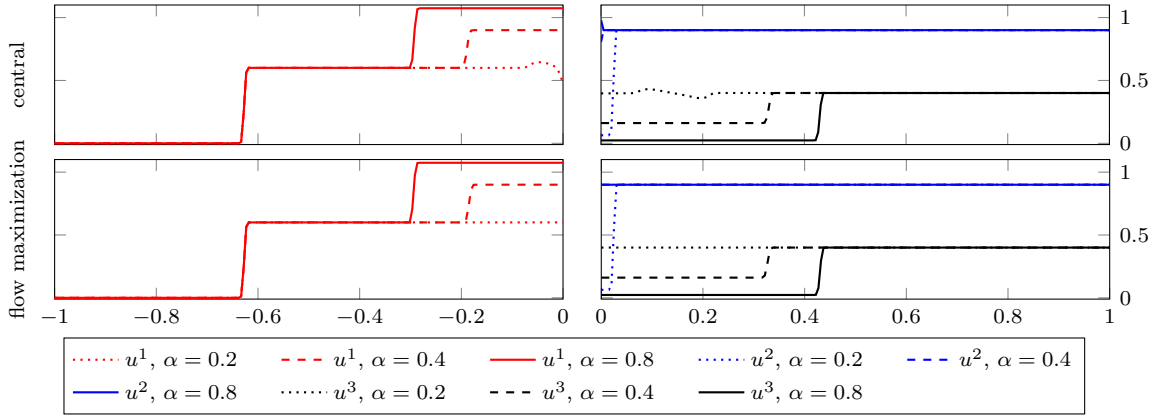

We employed the second order scheme on $200$ grid cells per edge with Courant number $\text{CFL=0.24}$ and relaxation speeds $\lambda_1=\lambda_2=\lambda_3=1$ to compute numerical solutions considering both the central relaxation and the flow maximization approach for varied $\alpha$.
Fig.~\ref{fig:LWR12} shows and compares these numerical solutions at time instance $T=0.75$. In case of larger $\alpha$ ($\alpha=0.4$ and $\alpha=0.8$ shown), the large vehicle density on road 2 leads to congestion and a backward moving traffic wave on the incoming road, whose amplitude increases as $\alpha$ increases. Moreover, the lack of incoming traffic leads to a shock wave on road 3. In case of small $\alpha$ ($\alpha=0.2$ shown) the traffic freely flows from the incoming to the outgoing roads. While the vehicle density on road 3 stays constant, a shock wave is formed on road 2 due to the lack of incoming traffic. For larger $\alpha$ the central approach reproduces the dynamics of flow maximization. In case of smaller $\alpha$ the central approach introduces small oscillations of low frequency in the numerical solutions of $u^1$ and $u^3$ close to the coupling node but still achieves qualitatively similar dynamics as the flow maximization approach.

\section{Conclusion}
The presented numerical experiments complement the numerical study in \cite{herty2022centr}. They show that the considered scheme allows for numerical simulation of networks with multiple outgoing edges even in case of non-convex flux functions. Further, the considered approach qualitatively reproduces the dynamics of flow maximization also in case of multiple outgoing edges. In addition, an experiment with the relaxation scheme and the LWR model indicates that the scheme is asymptotic preserving for some network problems.

\section*{Acknowledgements}
The author thanks the Deutsche Forschungsgemeinschaft (DFG, German Research Foundation) for the financial support through 320021702/GRK2326, 333849990/IRTG-2379, CRC1481, HE5386/18-1,19-2,22-1,23-1, ERS SFDdM035 and under Germany’s Excellence Strategy EXC-2023 Internet of Production 390621612 and under the Excellence Strategy of the Federal Government and the Länder. Support through the EU ITN DATAHYKING is also acknowledged.

\vspace{\baselineskip}








\providecommand{\WileyBibTextsc}{}
\let\textsc\WileyBibTextsc
\providecommand{\othercit}{}
\providecommand{\jr}[1]{#1}
\providecommand{\etal}{~et~al.}

\end{document}